\documentclass[12pt]{article}
\newtheorem{theorem}{Theorem}[section]
\newtheorem{cor}{Corollary}[section]
\newtheorem{prop}{Proposition}[section]
\newtheorem{lemma}{Lemma}[section]
\newcommand{\rf}[1]{\mbox{(\ref{#1})}}
\begin{document}

\title{Area minimization among marginally trapped surfaces in Lorentz-Minkowski space}
\author{By B{\footnotesize ENNETT} P{\footnotesize ALMER}}
\maketitle
\maketitle
\begin{abstract} We study an area minimization problem for spacelike zero mean curvature surfaces in four dimensional Lorentz-Minkowski space. The areas of these surfaces are compared of with the areas of certain marginally trapped surfaces having the same 
boundary values. 
\end{abstract}
\begin{center}
{\it
Department of Mathematics, Idaho State University,
Pocatello, ID 83209, U.S.A.
e-mail: palmbenn@isu.edu}
\end{center}
The quintessential  property of zero mean curvature surfaces in Euclidean space is that they locally minimize area with respect to their boundary curves. The method of calibrations can be used to show that any part of such a surface which can be represented as a graph over a convex domain, has least area when compared with any surface sharing the same boundary values. This justifies calling surfaces with zero mean curvature surfaces in Euclidean space `minimal surfaces'. In a similar way, zero mean curvature space-like surfaces in three dimensional Lorentz-Minkowski space are well known to locally maximize area.

 For space-like zero mean curvature surfaces in any four dimensional Lorentzian manifold there can be no analogous local minimizing or local  maximizing property. An arbitrarily  small neighborhood of any point possesses
an infinite dimensional space of deformations of the surface which fix the boundary and decrease the area and it has an infinite dimensional space of deformations which fix the boundary and increase the area. In a sufficiently small neighborhood of a point,
any compactly supported variation with space-like, (respectively timelike), variation field will initially increase, (respectively decrease), area. Space-like zero mean curvature surfaces are still characterized as equilibria for the area functional  so one would hope that they could be approached  in some meaningful way through variational methods.

We recall that a marginally trapped surface is a space-like surface whose mean curvature vector is isotropic (${\vec H}\cdot {\vec H}\equiv 0$).
It was shown in \cite{AP} that space-like zero mean curvature surfaces in four dimensional Lorentzian manifolds satisfying the Null Convergence Condition
 are weak local minima when compared with nearby marginally trapped surfaces having the same boundary values to first order. Specifically, each point has a neighborhood such that the second variation is non negative for all variations through marginally trapped surfaces having the same boundary values to first order. Consequently, a Morse index can be defined for any relatively compact subdomain of a space like zero mean curvature surface.   Marginally trapped surfaces were introduced in general relativity to study space-time singularities. They also occur naturally in the study of the conformal and Laguerre geometry
 of surfaces in Euclidean space.

In this paper we will study a new area minimizing property of zero-mean curvature surfaces in four dimensional Lorentz-Minkowski space which can roughly be described as follows.  Near   a `generic '  point on a space-like, zero mean curvature surfaces in ${\bf R}^4_1$, specifically a point where the Gauss curvature does not vanish, there are two  families of null direction fields of the normal space. We choose either one of them and consider the collection of marginally trapped surfaces having the same boundary values as $\Sigma$ and having mean curvature vectors lying in the chosen family of null directions. We show that $\Sigma$ has least area among this family of surfaces and we show that the family of surfaces is quite large.   The basic idea of the proof is to represent the surface using a gauge for which the equation for zero mean curvature is a linear elliptic fourth order equation. We represent the area as the Dirichlet energy for this equation and interpret the Dirichlet Principle for the energy as an area inequality.

It would be quite interesting to know if there is any analogue of this minimizing property for zero mean curvature surfaces in other four dimensional Lorentzian manifolds. The case where the ambient space is the deSitter space $S^4_1$ is of particular interest for two reasons.  
First of all, in this space area minimization
among marginally trapped surfaces has interesting applications in conformal geometry.
Secondly, the space ${\cal L}$ of null directions in $TS^4_1$  again has a product structure ${\cal L}\approx S^4_1\times S^2$ which might make the extension to this space accessible.

 This paper owes a great deal to the theory of Laguerre differential geometry as developed by Wilhelm Blaschke. For example, the representation formula and the connection between zero mean curvature surfaces in ${\bf R}^4_1$ and the equation \rf{DD}   appear in Blaschke's work \cite{B}.  In fact, our result can be interpreted as a local minimization property for Laguerre minimal surfaces. 
 
 This work was partially funded by Fellowship S-08154 from the Japan Society for the Promotion of Science.
 
 This work was also partially supported by Fundaci\'{o}n S\'{e}neca
project 04540/GERM/06, Spain. This research is a result of the activity
developed
within the framework of the Programme in Support of Excellence Groups of the
Regi\'{o}n de Murcia, Spain, by Fundaci\'{o}n S\'{e}neca, Regional
Agency for
Science and Technology (Regional Plan for Science and Technology 2007-2010).

\section{Main result}

We will call a space-like surface   $Y:\Sigma \rightarrow {\bf R}^4_1$ a {\it spherical graph} if there exists
a null section of the normal bundle $\perp Y$  of the form $\xi=(\nu, 1)$ such that $\nu:\Sigma \rightarrow S^2$ is injective.
If $Y$ is a spherical graph and if $\nu(\Sigma)=:\Omega$, we will say that $Y$ is a spherical graph over $\Omega$.
If a surface is a space-like spherical graph over $\Omega$, then the surface can be parameterized by $\nu^{-1}=:Y$ and
we write the surface as $Y:\Omega \rightarrow {\bf R}^4_1$.
If $Y:\Omega \rightarrow {\bf R}^4_1$ is a spherical graph over $\Omega$ and the mean curvature $ {\vec H}_Y$ of $Y$ satisfies
${\vec H}_Y\cdot (\nu, 1)\equiv 0$, then we will call $Y$ a marginally trapped spherical graph over $\Omega$. The reader is warned
that  we are not only requiring that the spherical graph is marginally trapped, but we are also prescribing the null direction of the mean curvature vector at points where it is non zero.
\begin{theorem}
\label{T1}
Consider a space-like zero mean curvature surface which can be represented as a spherical graph $X:\Omega  \rightarrow {\bf R}^4_1$.
Let $Y:\Omega \rightarrow {\bf R}^4_1$ be any marginally trapped spherical graph over $\Omega$ with
$Y|_{\partial \Omega}\equiv X|_{\partial \Omega}$ . Then,
$${\rm Area }[X]\le {\rm Area}[Y]\:,$$
holds.

   Furthermore,  if $\beta_1(\Omega)$ denotes the first  buckling eigenvalue of $\Omega$ (see \rf{buck} below for the definition),  then
   there holds
   \begin{equation}
   \label{est}
   (\beta_1(\Omega)-2)\int_\Omega |{\underline X}-{\underline Y}|^2\:d\omega
   +{\rm Area }[X]\le {\rm Area}[Y]\:.\end{equation}
   Here ${\underline X}$, (resp.$ {\underline Y}$) denotes the tangential part of the projection to ${\bf R}^3$ of $X$, (resp. $Y$).
\end{theorem}
{\it Remark}  It will be shown below that the problem of producing marginally trapped spherical graphs with prescribed boundary values is very underdetermined (Proposition \rf{p1}). Such surfaces can be produced from an arbitrary smooth function whose boundary values are prescribed to second order.\\[4mm]
{\it Remark} Note that any minimal surface in ${\bf R}^3\subset {\bf R}^4$ is a space like zero mean curvature surface. If the surface is oriented with normal $\nu$, then $(\nu, 1)$ is a null  section of the normal bundle of the surface. The usual catenoid can thus be considered as a zero mean curvature spherical graph since its Gauss map is globally injective. Sufficiently large pieces of the catenoid are unstable and hence they fail to minimize area with respect to their boundary. More precisely, for a sufficiently large domain bounded by two circles, there may exist a rescaling of the catenoid with the same boundary which has less area than the original piece.  This makes it clear that some additional hypothesis must be imposed on the comparison surfaces  to insure that a spherical graph with zero mean curvature is relatively area minimizing. Our definition achieves this by prescribing the Gaussian image of the surface.

We begin the proof with some preliminary results.

\begin{lemma}
\label{lemma1}
Let  $Y:\Sigma \rightarrow {\bf R}^4_1$ be a space-like surface. Assume that its mean curvature vector satisfies
${\vec H}\cdot (\nu, 1)\equiv 0$, with $|\nu|\equiv 1$ and $(\nu,1)\in \Gamma(\perp Y)$, Then the map  $\nu: \Sigma \rightarrow S^2$ is weakly conformal.
\end{lemma}
{\it Proof.} Let $\xi:=(\nu, 1)$. Then $\Delta Y\cdot \xi \equiv 0$. Since $Y_i\cdot \xi\equiv 0$, we have
$Y_{ii}\cdot \xi +Y_i\xi_i\equiv 0$. By choosing an orthonormal frame $\{e_i\}$ which diagonalizes $\nabla \xi$, we can write
$$\xi_i=\lambda_i Y_i +\alpha(e_i)\xi\:,$$
where $\alpha$ denotes the connection 1-form in the normal bundle $\perp Y$. Since $\lambda_1+\lambda_2=0$, we
have
\begin{equation}
\label{conf}
d\xi \cdot d\xi =\lambda_1^2\: dY\cdot dY\:.\end{equation}
{\bf q.e.d.}\\[4mm]
We will now consider a spherical graph $\Sigma \rightarrow {\bf R}^4_1$  whose mean curvature satisfies ${\vec H}\cdot (\nu,1)\equiv 0$.
 By the lemma, we
can assume that the surface is represented as a conformal immersion
$Y:=\nu^{-1}:\Omega\subset S^2\rightarrow  {\bf R}^4_1$. Again, let $\xi=(\nu, 1)$ and define
$f:=Y\cdot \xi$.

We denote the Laplacian on $S^2$ by ${\hat \Delta}$. The gradient of a function $u$ on $S^2$ will be denoted $Du$.
\begin{lemma}
\begin{equation}
\label{rep}
Y=(Df-(\frac{1}{2}{\hat \Delta}f)\nu \:,\:-\frac{1}{2}[{\hat \Delta}f+2f])\:.\end{equation}
\end{lemma}
{\it Proof.} We will write a vector in ${\bf R}^4_1$ as $v=({\underline v}, v_4)$.
Let $e_r$ be an orthonormal frame defined on a neighborhood in $S^2$.
We have $f_r=Y\cdot \xi_r={\underline Y}\cdot e_r\:$. This shows that modulo $\nu$,
${\underline Y}=Df$ holds. Now calculate
$$f_{rr}=Y\cdot \xi_{rr}+Y_r\cdot \xi_r= {\underline Y}\cdot \nu_{rr}+Y_r\cdot \xi_r\:.$$
Sum this over $r$ using that ${\hat \Delta }\nu=-2\nu$ and that $\sum Y_r\cdot \xi_r=0$
since the mean curvature vector of $Y$ is parallel to the null direction $\xi$. We obtain
\begin{equation}
\label{eq1}
{\hat \Delta }f=-2Y\cdot (\nu, 0)=-2Y\cdot (\nu,1)+2Y\cdot ({\underline 0}, 1)\:.\end{equation}
Thus
\begin{equation}
\label{eq2}  ({\hat \Delta}+2)[f]=2Y\cdot ({\underline 0}, 1)\:\end{equation}
holds. Combining this with the first equality in \rf{eq1}, gives \rf{rep}. {\bf q.e.d}
\begin{cor}
$Y$ has zero mean curvature only if ${\hat \Delta }({\hat \Delta }f+2f)=0$ holds.
\end{cor}
The necessity is clear from \rf{eq2} since ${\vec H}_Y=0$ if and only if ${\hat \Delta Y}=0$.
\begin{lemma} Let $Y$ be as above. Then the area of the surface is given by
\begin{equation}
\label{area}
{\rm Area}\:[Y]=\frac{1}{4}\int_\Omega {\hat \Delta}[f]\:( {\hat \Delta}+2)[f]\:d\omega +
\frac{1}{2}\oint_{\partial \Omega} \frac{1}{2} \partial_n|Df|^2-( {\hat \Delta}+2)[f]\partial_nf\:d\sigma\:.\end{equation}
\end{lemma}
{\it Proof.} We can rewrite \rf{rep} as
\begin{equation}
\label{alt}
Y=(Df +f\nu, 0)-\frac{1}{2}[{\hat \Delta}f+2f])(\nu, 1)\:.
\end{equation}
If we let $I=d\nu$ denote the identity endomorphism on each tangent space, we have
$$dY=(D^2f+fI, 0)-\frac{1}{2}d\bigl(( {\hat \Delta}+2)[f]\bigr)\otimes(\nu,1)-\frac{1}{2}\bigl(( {\hat \Delta}+2)[f]\bigr)(I,0)\:.$$

Rearranging terms, we obtain

\begin{equation}
\label{dY}
dY=(D^2f-\frac{1}{2}{\hat \Delta}[f]\otimes I-\frac{1}{2}d\bigl(( {\hat \Delta}+2)[f]\bigr)\otimes \nu,
-\frac{1}{2}d\bigl(( {\hat \Delta}+2)[f]\bigr)\:.\end{equation}
Here $I$ denotes the identity endomorphism field on $TS^2$.
From this it follows that the metric induced by $Y$ is given by
\begin{equation}
\label{g}
dY\cdot dY=(D^2f-\frac{1}{2}{\hat \Delta}[f]\otimes I)\cdot (D^2f-\frac{1}{2}{\hat \Delta}[f]\otimes I)\:.\end{equation}
Note that $(D^2f-\frac{1}{2}{\hat \Delta}[f]\otimes I)$ is the trace free part of the Hessian
$D^2f$. Recalling the formula \rf{conf}, we have that the metric appearing in \rf{g} is
conformal to the metric $d{\cal S}^2$  on $S^2$.

Let $\sigma_i$ be the eigenvalues of $D^2f$. Then from \rf{g}, we have
\begin{eqnarray*}
dY\cdot dY&=&\frac{1}{4}(\sigma_1-\sigma_2)^2d{\cal S}^2\\
                 &=&\frac{1}{4}\bigl((\sigma_1+\sigma_2)^2-4\sigma_1\sigma_2\bigr)d{\cal S}^2\\
                 &=&\frac{1}{4}\bigl( ({\hat \Delta}[f])^2-4M[f]\bigr)d{\cal S}^2\:,
 \end{eqnarray*}
 where $M[f]$ is the Monge-Ampere operator (determinant of the Hessian).  It then follows that
$${\rm Area}\:[Y]=\frac{1}{4}\int_\Omega ( {\hat \Delta}[f] )^2-4M[f]\:d\omega\:.$$

We now recall a version of an integrated Lichnerowicz formula. If $\Sigma$ is any surface with boundary and
$u$ is any sufficiently smooth function
$$\oint_{\partial \Sigma} \frac{1}{2}\partial_n |\nabla u|^2-(\Delta u)\partial_n u \:ds=
\int_\Sigma |\nabla^2 u|^2  -(\Delta u)^2+K|\nabla u|^2\:d\Sigma\:.$$
Note that $ |\nabla^2 u|^2  -(\Delta u)^2=-2M[u]$.
We apply this to $f$ using our current notation and using $K\equiv 1$ on $S^2$,  to obtain
$${\rm Area}\:[Y]=\frac{1}{4}\int_\Omega ( {\hat \Delta}[f] )^2-2|Df|^2\:d\omega+\frac{1}{2}\oint_{\partial \Omega} \frac{1}{2}\partial_n|Df|^2-{\hat \Delta}[f]\partial_nf\:ds\:.  $$
Finally, using
$$\int_\Omega |Df|^2\:d\omega =\oint_{\partial \Omega} f\partial_nf\:ds -\int_\Omega f{\hat \Delta}[f]\:d\omega\:,$$
we obtain \rf{area}.{\bf q.e.d.}
\begin{prop}
\label{p1}
Let $f$ be any smooth function on a domain $\Omega \subset S^2$ and define $Y=Y_f:\Omega \rightarrow {\bf R}^4_1$ by
\rf{rep}. Then $Y$ is marginally trapped. Consequently, if $g$ is any smooth function on $\Omega$ such that $f$ and $g$ agree to second order on
$\partial \Omega$, then $Y_g$ is a marginally trapped spherical graph over $\Omega$ which agrees with $Y_f$  on $\partial \Omega$.
\end{prop}
{\it Proof.}  From \rf{dY}, it follows that
$${\rm trace}\:( dY\cdot d(\nu,1))={\rm trace} \:((D^2f-\frac{1}{2}{\hat \Delta}[f]\otimes I)\cdot I )=0\:.$$
This shows that $Y$ is marginally trapped.

The second statement follows immediately from formulas \rf{rep} and \rf{dY}. {\bf q.e.d.}

Consider the expression for the marginally trapped surfaces given by \rf{alt}.
The equation ${\bar X}:=Df+f\nu:\Omega \subset S^2\rightarrow {\bf R}^3$ represents 
a surface in three dimensional space having Gauss map $\nu$. The curvature of this surface is never zero since its Gauss map is just $X^{-1}$.  Also ${\hat \Delta f}+2f=
-(k_1^{-1}+k_2^{-1})$, where $k_i$ are the principal curvatures if $X$. Thus, the immersion ${\bar X}$ encodes all of the extrinsic geometry of the map $Y$ and gives us a way to visualize the marginally trapped surface. 
%\begin{figure}
   %\includegraphics[width=80mm,height=80mm,angle=0]{slpic.jpg}
     % \caption{Small liquid drop}
% \end{figure}
\begin{lemma}
\label{be}
For $\Omega \subset\subset S^2$, there  holds.
\begin{equation}
\label{buck}\inf_{u\in W^{2,2}_0(\Omega)\atop{ u\neq const}}\frac{\int({\hat \Delta}u)^2\:d\omega}{\int |Du|^2\:d\omega}=:\beta_1(\Omega) \ge 2\:.
\end{equation}
\end{lemma}

The number $\beta_1(\Omega)$ appearing above is called the first buckling eigenvalue.
The space $ W^{2,2}_0(\Omega)$ is the completion of $C^\infty_C(\Omega)$ with respect to the metric
$$||f||_{ W^{2,2}} =\int_\Omega f^2+|Df|^2+|D^2f|^2\:d\omega\:.$$
A function $f \in C^\infty$ is in $W^{2,2}_0{\Omega}$ if and only if $f\equiv 0$ and 
$Df\equiv 0$ on $\partial \Omega$. 

{\it Proof.} It is easy to see that the infimum of the same quotient over all smooth non constant functions on $S^2$ is
2. This can be seen by expanding $u$ as a series of eigenfunctions of the Laplacian. If there is a subdomain
$\Omega$ and a $W^{2,2}_0$ function which makes the quotient negative, then extending this function to be identically
zero off $u$ leads to a contradiction. {\bf q.e.d.}
\begin{lemma}(Dirichlet's Principle). If $q$ solves
\begin{equation}
\label{EL}
{\hat \Delta }({\hat \Delta }q+2q)=0\end{equation}
 in $\Omega \subset \subset S^2$ and $f$ is any smooth function with $f-q\in W^{2,2}_0({\Omega})$, then
\begin{equation}
\label{dirichlet}
(\beta_1(\Omega)-2)\int_\Omega |Df-Dq|^2\:d\omega+\int_\Omega {\hat \Delta}[q]\:( {\hat \Delta}+2)[q]\:d\omega\le \int_\Omega {\hat \Delta}[f]\:( {\hat \Delta}+2)[f]\:d\omega\:,\end{equation}
holds.\end{lemma}
{\it Proof.} After a partial integration, the previous lemma can be interpreted as saying that for any non constant $W_0^{2,2}$
on $S^2$,
$$\int {\hat \Delta}[u]\cdot ({\hat \Delta}+2)[u]\:d\omega \ge (\beta_1(\Omega)-2)\int_\Omega |Du|^2\:d\omega\:,$$
holds.

Since $f-q\in W^{2,2}_0({\Omega})$ holds, we can apply Lemma \rf{be}, to obtain
\begin{eqnarray*}
(\beta_1-2)\int_\Omega |Df-Dq|^2\:d\omega&\le& \int_\Omega {\hat \Delta}[f -q]\cdot ({\hat \Delta}+2)[f-q]\:d\omega\\
  &=&\int_\Omega {\hat \Delta}[f] \cdot ({\hat \Delta}+2)[f]\:d\omega+\int_\Omega {\hat \Delta}[q] \cdot ({\hat \Delta}+2)[q]\:d\omega\\
  &&-\int_\Omega {\hat \Delta}[q] \cdot ({\hat \Delta}+2)[f]\:d\omega-\int_\Omega {\hat \Delta}[f] \cdot ({\hat \Delta}+2)[q]\:d\omega\\
  &=&(*)
  \end{eqnarray*}
The condition  $f-q\in W^{2,2}_0({\Omega})$ means that $f\equiv q$ and $\partial_nf\equiv \partial_nq$
on $\partial \Omega$. By using Green's second identity and \rf{EL}, we get
 \begin{eqnarray*}
 \int_\Omega {\hat \Delta}[q] \cdot ({\hat \Delta}+2)[f]\:d\omega&=&
 \int_\Omega {\hat \Delta}[q] \cdot ({\hat \Delta}+2)[f]\:d\omega-\int_\Omega f\cdot {\hat \Delta}({\hat \Delta}+2)[q]\:d\omega\\
 &=& \int_\Omega {\hat \Delta}[q] \cdot {\hat \Delta}[f]-f{\hat \Delta}^2[q]\:d\omega\\
 &=&\oint_{\partial \Omega} {\hat \Delta}[q]\partial_n f -f\partial_n{\hat \Delta}[q]\:ds\\
  &=&\oint_{\partial \Omega} {\hat \Delta}[q]\partial_n q -q\partial_n{\hat \Delta}[q]\:ds\\
   &=& \int_\Omega {\hat \Delta}[q] \cdot {\hat \Delta}[q]-q{\hat \Delta}^2[q]\:d\omega\\
   &=& \int_\Omega {\hat \Delta}[q] \cdot ({\hat \Delta}+2)[q]\:d\omega\\
   &=&(**)
 \end{eqnarray*}
 Next note that
\begin{eqnarray*}
 \int_\Omega {\hat \Delta}[f] \cdot ({\hat \Delta}+2)[q]\:d\omega&=& \int_\Omega {\hat \Delta}[f] \cdot ({\hat \Delta}+2)[q]\:d\omega\\
                                                                                                   &&+2\int_\Omega f{\hat \Delta}[q]-q{\hat \Delta}[f]\:d\omega\\
                                                                                                   &=&(***)
 \end{eqnarray*}
 The last integral on the right is zero since
$$\int_\Omega q{\hat \Delta}[f]-f{\hat \Delta}[q]\:d\omega=\oint_{\partial \Omega} q\partial_nf-f\partial_nq\:ds =
\oint_{\partial \Omega} q\partial_nq-q\partial_nq\:ds\:.$$
Combining (*), (**) and (***), we get
$$(\beta_1(\Omega)-2)\int_\Omega |Df-Dq|^2\:d\omega\le \int_\Omega {\hat \Delta}[f] \cdot ({\hat \Delta}+2)[f]\:d\omega-\int_\Omega {\hat \Delta}[q] \cdot ({\hat \Delta}+2)[q]\:d\omega\:,$$
which is \rf{dirichlet}.\\[4mm]
{\it Proof of Theorem \rf{T1}}
We can assume that $X$ is a zero mean curvature spherical graph which is represented over $\Omega \subset S^2$ using formula \rf{rep} with $f$ replaced by a function $q$ satisfying equation \rf{EL}. If $Y$ is
a marginally trapped spherical graph over $\Omega$ with $Y\equiv X$ on the boundary we also represent $Y$
using \rf{rep} with $f=Y\cdot(\nu,1)$. We claim that $f-q\in W^{2,2}_0(\Omega)$.

It is clear that $f\equiv q$ on $\partial \Omega$ since $q=X\cdot (\nu,1)$.  By projecting both $X$ and $Y$
to the three dimensional space and using that $\nu$ is perpendicular to both $Df$ and $Dq$, we see that
$Df\equiv Dq$ along $\partial \Omega$ also, proving the claim.

It suffices to show that 
if $q$ solves \rf{EL} and $f-q\in W^{2,2}_0{\Omega}$, then
\begin{equation}
\label{area2}
{\rm Area}\:[X]\le {\rm Area}\:[Y]\:,
\end{equation}
where
$$X=(Dq-(\frac{1}{2}{\hat \Delta}q)\nu \:,\:-(\frac{1}{2}[{\hat \Delta}q+2q]),\qquad Y=(Df-(\frac{1}{2}{\hat \Delta}f)\nu \:,\:-\frac{1}{2}[{\hat \Delta}f+2f]).$$

By the previous lemma, the result will follow from \rf{area} if we can show that the boundary integrals in \rf{area} is unchanged
if we replace $f$ by $q$.  We denote the unit tangent and normal to $\partial \Omega$ by $t$ and $n$.  We  have
\begin{eqnarray*}
 \frac{1}{2} \partial_n|Df|^2-( {\hat \Delta}+2)[f]\partial_nf&=&[D_n(Df)\cdot Df]-(D_tDf\cdot t+D_nDf\cdot n)f_n-2ff_n\\
&=&[(D_n(Df)\cdot t)f_t+(D_n(Df)\cdot n)f_n]-(D_tDf\cdot t+D_nDf\cdot n)f_n\\
&&-2ff_n\\
&=&(D_n(Df)\cdot t)f_t-(D_tDf\cdot t)f_t-2ff_n\\
&=&(D_t(Df)\cdot n)f_t-(D_tDf\cdot t)f_t-2ff_n,
 \end{eqnarray*}
 where we have used that the Hessian is symmetric. By assumption $f \equiv q$, $Df \equiv Dq$ on $\partial \Omega$, and thus
 $f_n \equiv q_n$ and $D_tDf \equiv D_tDq$ on $\partial \Omega$. This gives the result. {\bf q.e.d.}
 \bigskip

 We now want to relate the local geometry of a surface with ${\vec H}=0$ to the property of being a spherical graph. We work locally and introduce
 a complex coordinate on a neighborhood of $\Sigma$. We write the induced metric as $ds_X^2=e^\rho |dz|^2$. Let $\{\xi, \eta\}$  be a local frame for the normal bundle with $\xi \cdot \xi \equiv 0\equiv \eta \cdot \eta$ and $\xi \cdot \eta \equiv 1$. We can express the usual Frenet equations of the immersion as
 \begin{eqnarray*}
 X_{zz}&=& \rho_z  X_z+ \frac{\phi}{2} \xi+ \frac{\psi}{2}\eta\\
 X_{z{\bar z}}&=&0\\
 \xi_z&=& -{\psi}e^{-\rho}X_{\bar z}+\sigma \xi \\
 \eta_z&=& -{\phi}e^{-\rho}X_{\bar z}-\sigma \eta \:,
 \end{eqnarray*}
 where $\phi dz^2$ and $\psi dz^2$ are invariantly defined quadratic differentials.
 The Gauss equation becomes
\begin{equation}
\label{curv}
K=-e^{-2\rho}(\phi {\bar \psi}+{\bar \phi}\psi)\:,\end{equation}
and the Codazzi equations are
\begin{equation}
\phi_{\bar z}=-{\bar \sigma}\phi\:\qquad \psi_{\bar z}={\bar \sigma }\psi \:.
\end{equation}
It follows that both $\phi$ and $\psi$ are pseudoholomorphic. 

Zero mean curvature surfaces in ${\bf R}^4_1$ can be parameterized using a Weierstrass representation which depends on holomorphic functions. It follows from this that these surfaces are real analytic and hence the curvature is either identically zero or the points where $K=0$ are nowhere dense.

For some functions $\alpha$ and $\beta$, we can write the null sections as
\begin{equation}
\label{x}
\xi=\alpha(\nu_1,1)\:,\qquad \eta=\beta(\nu_2,1)\:,
\end{equation}
where $\nu_i$ are maps into $S^2$. We have
$$(d\nu_1,0)^T= \alpha \nabla^T(\nu_1,,1)=  \nabla^T\xi = -2\Re(\phi e^{-\rho} X_{\bar z}\otimes dz)\:,$$
$$(d\nu_2,0)^T= \beta \nabla^T(\nu_2,,1)=  \nabla^T\eta = -2\Re(\psi e^{-\rho} X_{\bar z}\otimes dz)\:.$$
Thus is $d\nu_i(p)=0$ at some point $p$, then either $\phi(p)=0$ or $\psi(p)=0$ and by \rf{curv}, $K(p)=0$ also.
It follows that if $K(p)\ne 0$, then the surface is, locally,  a spherical graph over domains in $S^2$ defined by
the maps $\nu_i$ coming from both null directions in the normal bundle.

We end with a remark about solving the equation
\begin{equation}
\label{DD}
{\hat \Delta }({\hat \Delta }f+2f)=0\:.\end{equation}
These remarks can be found in \cite{B}.
Let $\Omega$ be a domain in $S^2$ and write the usual metric on $S^2$ as $d{\cal S}^2=e^\sigma |dw|^2$ where
$w$ is a complex coordinate. The Gauss equation gives:
$$\sigma_{w{\bar w}}=\frac{-1}{2}e^\sigma\:.$$
Let $f$ be a holomorphic function on $\Omega$, i.e. $f_{\bar w}=0$. It is straightforward to derive that that the function
$Q:=(-1/2)\sigma_w f$ satisfies 
$${\hat \Delta Q}=4e^{-\sigma}Q_{{w}{\bar w}}=f_w+\sigma_w f\:.$$
It is also straightforward to derive that 
$$({\hat \Delta}+2)(f_w+\sigma_wf)=0\:,$$
whenever $f$ is holomorphic. It follows that $Q$ is a complex solution of \rf{DD}. Since the operator 
${\hat \Delta }({\hat \Delta }+2)$ is real, the real and imaginary parts  $\Re(Q)$ and $\Im(Q)$ are real solutions of \rf{DD}.
Evidently, all zero mean curvature surfaces produced from these solutions have non vanishing curvature and are hence locally area minimizing with respect to marginally trapped surfaces having the same first order boundary values.


\begin{thebibliography}{2}
\bibitem{AP} Al\'{\i}as, L. J. and Palmer, B.,  {\it Deformations of Stationary Surfaces}.
  Classical and Quantum Gravity, 14 (1997), 2107--2111.
\bibitem{B} Blaschke, W., Laguerre geometrie III, Beitrage zur Flachentheorie,
Hambg. Abh. 4(1926), 1-12.
\bibitem{P} Pottman, H. , Grohs, P, Mitra, N. J. , {\it Laguerre  minimal surfaces; isotropic geometry and linear elasticity}. To appear in Advances in Computational Mathematics.
\end{thebibliography}
\end{document}